
\tolerance=10000
\raggedbottom

\baselineskip=15pt
\parskip=1\jot

\def\sk{\vskip 3\jot}

\def\heading#1{\vskip3\jot{\noindent\bf #1}}
\def\label#1{{\noindent\it #1}}


\def\ref#1;#2;#3;#4;#5.{\item{[#1]} #2,#3,{\it #4},#5.}
\def\refinbook#1;#2;#3;#4;#5;#6.{\item{[#1]} #2, #3, #4, {\it #5},#6.} 
\def\refbook#1;#2;#3;#4.{\item{[#1]} #2,{\it #3},#4.}


\def\({\bigl(}
\def\){\bigr)}


\def\ga{\gamma}
\def\ph{\phi}
\def\me{\omega}

\def\Si{\Sigma}
\def\Ph{\Phi}
\def\Ps{\Psi}


\def\nunit{\overline{1}}

\def\Ex{{\rm Ex}}
\def\Var{{\rm Var}}

\def\parity{{\rm parity}}
\def\majority{{\rm majority}}

\def\rmi{{\rm i}}
\def\rmii{{\rm ii}}

\def\abs#1{\vert#1\vert}

\def\bfC{{\bf C}}

{
\pageno=0
\nopagenumbers

\rightline{\tt izsak.tex}
\vskip1in

\centerline{\bf Carry Propagation in  Multiplication by Constants}
\vskip0.5in

\centerline{Alice Izsak}
\centerline{\tt aizsak@gmail.com}
\vskip0.15in

\centerline{Department of Computer Science}
\centerline{University of British Columbia}
\centerline{2366 Main Mall}
\centerline{Vancouver, BC V6T 1Z4}
\centerline{CANADA}
\vskip0.5in

\centerline{Nicholas Pippenger}
\centerline{\tt njp@math.hmc.edu}
\vskip0.15in

\centerline{Department of Mathematics}
\centerline{Harvey Mudd College}
\centerline{1250 Dartmouth Avenue}
\centerline{Claremont, CA 91711}
\vskip0.5in

\noindent{\bf Abstract:}
Suppose that a random $n$-bit number $V$ is multiplied by an odd constant $M\ge 3$, by adding shifted versions of the number $V$ corresponding to the $1$s in the binary representation of the constant $M$.
Suppose further that the additions are performed by carry-save adders until the number of summands
is reduced to two, at which time the final addition is performed by a carry-propagate adder.
We show that in this situation the distribution of the length of the longest carry-propagation chain in the final addition is the same (up to terms tending to $0$ as $n\to \infty$) as when two independent $n$-bit numbers are added,
and in particular the mean and variance are the same (again up to terms tending to $0$).
This result applies to all possible orders of performing the carry-save additions.
It also applies if the constant multiplier is recoded to reduce the number 
of operations by allowing subtractions as well as additions, corresponding to occurrences of 
a ``negative unit'' $\nunit  = -1$ in a representation of $M$.
\vfill\eject
}

\heading{1.  Introduction}

Let $X$ and $Y$ be random $n$-bit integers that are independent and uniformly distributed
in $[0,2^n - 1]$.
If they are added in the usual way, starting at their rightmost end and proceeding to the left,
their may be various ``carry-propagation chains''.
A {\it carry-propagation chain\/} is a sequence of $k\ge 1$ consecutive positions in the binary
representations of $X$ and $Y$ in which the rightmost position {\it generates\/} a carry
(because both $X$ and $Y$ contain $1$s in these positions), and the remaining $k-1$
positions to the left {\it propagate\/} this carry (because one, but not the other, of $X$ and $Y$ contains
a $1$ in each of these positions).
Let the random variable $\bfC_n$ denote the length of the longest carry-propagation chain.
(Note that the longest carry-propagation chain is not necessarily the longest sequence of consecutive carries: the addition of the binary numbers $0101$ and $1111$ gives rise to two carry-propagation chains, each of length two, not to one of length four.)
The length of the longest carry-propagation chain is of interest because it governs the execution
of certain {\it parallel\/} implementations of addition (see Claus [C] and Knuth [K]).

The distribution of $\bfC_n$ has been investigated since the early days of electronic computing.
The investigation was begun in the famous report of Burks, Goldstein and von Neumann
[B] in 1946, where it was shown that $\Ex(\bfC_n) \le \log_2 n + 1$.
The next step was taken by Claus [C], who showed that $\Ex(\bfC_n)\ge \log_2 n - 2$.
Knuth [K] showed that 
$$\Pr(\bfC_n\ge k) = 1 - e^{-n/2^{k+1}} + O\left({(\log n)^3\over n}\right) \eqno(1.1)$$
(where the constant in the $O$-term is independent of $k$ as well as $n$),
and that this implies
$$\Ex(\bfC_n) = \log_2 n +  \ga\log_2 e - {3\over 2} - \Ph(\log_2 n) + 
O\left({(\log n)^4\over n}\right), \eqno(1.2)$$
where $\ga = 0,5772\ldots$ is Euler's constant, $e = 2.718\ldots$ is the base of natural logarithms,
and $\Ph(\nu)$ is a periodic function of $\nu$ with period $1$ and average $0$
(that is, $\int_0^1 \Ph(\nu)\,d\nu = 0$) satisfying 
$\abs{\Ph(\nu)}\le 1.573\ldots\times 10^{-6}$ for all $\nu\in[0,1)$.
Pippenger [P] gave an elementary derivation of (1.1), and showed that it also implies
$$\Var(\bfC_n) = {\pi^2\over 6}(\log_2 e)^2 + {1\over 12} + \me + \Ps(\nu) + 
O\left({(\log n)^5\over n}\right),\eqno(1.3)$$
where $\pi = 3.14159\ldots$ is the circular ratio, $\me = 1.2374\ldots\times 10^{-12}$ is a constant,
and $\Ps(\nu)$ is a periodic function of $\nu$ with period $1$ and average $0$ satisfying 
$\abs{\Ps(\nu)}\le 5.3573\ldots\times 10^{-6}$ for all $\nu\in[0,1)$.

In Section 2 we shall present a new analysis of the addition problem that yields results similar to
those above, but with weaker error bounds.
Specifically, we shall show that
$$\Pr(\bfC_n\ge k) = 1 - e^{-n/2^{k+1}} + O\left({\log n\over n^{1/3}}\right). \eqno(1.4)$$
This implies
$$\Ex(\bfC_n) = \log_2 n +  \ga\log_2 e - {3\over 2} - \Ph(\log_2 n) + 
O\left({(\log n)^2\over n^{1/3}}\right) \eqno(1.5)$$
and
$$\Var(\bfC_n) = {\pi^2\over 6}(\log_2 e)^2 + {1\over 12} + \me + \Ps(\nu) + 
O\left({(\log n)^3\over n^{1/3}}\right) \eqno(1.6)$$
in the same way that (1.1) implies (1.2) and (1.3).
The weaker error bounds are a result of our choice to present our new argument in its simplest form;
these bounds could be improved by elaboration of the argument (but, as Knuth [K] points out, so could those of (1.1--3)).
Our motivation, however, for presenting this new analysis is that it  can be extended to obtain the results
claimed in the abstract, which we shall now describe in more detail.

We shall investigate the length of the longest carry propagation chain that occurs when a random 
$n$-bit integer $V$, uniformly distributed in $[0, 2^n - 1]$, is multiplied by a fixed constant $M$.
The simplest case of our problem is $M=3$.
In this case, the product $Z = M\cdot V$ is obtained by adding $V$ to the number $2V$ that is obtained by shifting $V$ one position to the left.
The two random numbers being added in this case are {\it not\/} independent, but Izsak [I] has shown 
that the length of the longest carry-propagation chain nevertheless satisfies the estimate (1.1).
More generally, we may consider the case $M = 2^d + 1$ (where $d\ge 1$), for which
the product $Z = M\cdot V$ is obtained by adding $V$ to the number $2^d \,V$ that is obtained by shifting $V$  to the left $d$ positions.
Izsak [I] has shown that again the estimate (1.1) applies (where now the constant in the $O$-term may depend on $d$, but not on $k$ or $n$).

We shall consider a further generalization in which $M$ has two {\it or more\/} $1$s in its binary representation.
Suppose that the binary representation of $M$ is $M = \sum_{0\le j\le d} m_j \, 2^j$ (with
$m_j \in\{0,1\}$) and that $c$ (where $2\le c\le d+1$) of the digits $m_0, m_1, \ldots, m_d$ are $1$s
(so that the remaining $d+1-c$ are $0$s).
We may assume without loss of generality that $m_d = 1$ (since otherwise we could reduce the value of 
$d$) and that $m_0 = 1$ (since the carries that occur when multiplying by $2M$ will just be shifted versions of those that occur when multiplying by $M$).
Let $s_1 = 0 < s_2 < \cdots < s_c = d$ be the positions of the $1$-bits, so
$M = \sum_{1\le i\le c} 2^{s_i}$.
For $1\le i\le c$, let $W_i = 2^{s_i}\,V$ be obtained by shifting $V$ to the left $s_i$ positions. 
The product $Z = M\cdot V$ will be obtained by adding these $c$ numbers:
$Z = \sum_{1\le i\le c} W_i$. 

When $c = 3$, we can form the sum $Z = W_1 + W_2 + W_3$ in two {\it stages\/} as follows.
The first stage will perform a ``carry-save addition'', which takes the three numbers $W_1$, $W_2$ and $W_3$ and inputs and produces as outputs two numbers $X$ and $Y$ having the same sum:
$X + Y =  W_1 + W_2 + W_3 $.
There are of course many pairs of numbers $X$ and $Y$ that satisfy this condition.
The details of carry-save addition, including the specification of the numbers $X$ and $Y$ that will be produced, will be given later.
For now we merely observe that
in carry-save addition, all carries propagate one position to the left, and in a parallel implementation, all carries propagate simultaneously, so that a carry-save addition contributes a {\it fixed\/} delay to the parallel execution time.
Thus our analysis will not deal with carries in this stage.
The second stage will perform a conventional ``carry-propagate addition'' to obtain the final product $Z$ as the sum of $X$ and $Y$.
This addition is analogous to those considered in previous paragraphs, and it is the carry-propagation chains in this stage that will be the focus of our analysis.
We will obtain the estimate (1.4).

When $c\ge 4$, we can use $c-2$ carry-save additions to reduce the $c$ numbers
$W_1, W_2, \ldots W_c$ to two numbers $X$ and $Y$ in the first stage,
then add these two numbers with a carry-propagate addition
in the second stage to obtain $Z$ as before.
In this case, however, there is an additional complication:
there is more than one way to use $c-2$ carry-save additions to reduce $c$ numbers to two numbers.
At one extreme, one can sum $W_1$, $W_2$ and $W_3$ with the first cary-save addition,
then proceed similarly with the resulting $(c-3) + 2 = c-1$ numbers, and so forth.
The numbers $X$ and $Y$ are thus obtained after $c-2$ carry-save additions, each (except for the first)of which depends for at least one of its inputs on the its predecessor, so that these carry-save additions
contribute $c-2$ fixed delays to the parallel execution time.
At the other extreme, one can use $\lfloor c/3\rfloor$ carry-save additions in parallel
to combine $3\,\lfloor c/3\rfloor$
numbers, producing $2\,\lfloor c/3\rfloor$ numbers having the same sum, then proceed similarly with the resulting $(c - 3\,\lfloor c/3\rfloor) + 2\,\lfloor c/3\rfloor = c - \lfloor c/3\rfloor$ numbers, and so forth.
As Wallace [W] has observed, these $c-2$ carry-save additions
contribute only $\log_{3/2} c + O(1)$ fixed delays to the parallel execution time.
Our result, which is that the estimate (1.4) again holds for the carry-propagate addition in the second stage, applies equally to all of the ways of performing the carry-save addition in the first stage.

Finally, we shall consider a further generalization in which the number of carry-save additions
is reduced by using subtractions as well as additions.
For example, if $M=7$, we can represent $M$ in ``extended binary'' as $100\nunit$ rather than $111$,
where the ``negative unit'' $\nunit = -1$  means that the corresponding power of $2$ should be 
{\it subtracted\/} rather than added.
The product  $Z = 7V$ can then be computed by immediately performing a carry-propagate addition of 
$8V$ and $-V$, rather than by first performing a carry-save addition to combine $4V$, $2V$ and $V$,
and then performing a carry-propagate addition on the results.
(As this example shows, it may be necessary to increase $d$ by one when using an extended
binary representation.
But this merely increases the largest shift length, which does not affect the parallel execution time.
The goal is to reduce the number $c$ of summands, which {\it can\/} affect the parallel execution time.)
In general, we will want to find the extended binary representation of $M$
(that is, $M = \sum_{0\le i\le d} m_i \, 2^i$ with $m_i \in\{-1,0,1\}$, $m_d = 1$ and $m_0 = \pm 1$)
that minimizes the the number of non-zero digits $m_i = \pm 1$.
One representation that accomplishes this minimization is the ``canonical'' representation,
described by Lehman [L1, L2], Tocher [T] and Reitwiesner [R].
But non-canonical representations may tie the canonical one in achieving the minimum
(for example, the conventional $11$ ties the canonical $10\nunit$ in representing $3$),
so there may again be more than one optimal representation.
Our result, the estimate (1.4), again holds for all extended representations, optimal and non-optimal.

All of our results reinforce one point: the randomness in one uniformly distributed number $V$
is sufficient to produce the distribution (1.4); the full power of the independence of $X$ and $Y$ in 
forming their sum is not needed.
In Section 3, we shall give a specification at the bit level of the algorithms that were specified
above at the level of operations on numbers, and describe the features, common to all these algorithms, that will be used in the subsequent analysis.
In Section 4, we shall give the proof of (1.4) based, on these common features.
\sk

\heading{2.  A New Analysis of Addition}

In this section, we shall prove (1.4) for the addition of two independent random numbers.
The analyses of Knuth [K] and Pippenger [P] of (1.1) proceed by deriving a recurrence for the probability
that the addition of two random $n$-bit numbers yields a carry propagation chain of length at least $k$,
then solving this recurrence for the asymptotic behavior of this probability.
Our new analysis is based on the observation that the main term $1 - e^{-n/2^{k+1}}$ in (1.1) and (1.4)
is the probability that a Poisson-distributed random variable with mean $n/2^{k+1}$ has value at least one.
There are approximately $n$ (actually $n - k + 1$) places at which a carry-propagation chain of length $k$ can occur, and the probability that such a chain occurs at a given place is $1/2^{k+1}$.
If all these possible occurrences were independent, we could derive the desired result from the Poisson approximation to binomial distribution.
They are not independent, but the effects of their dependence can be analyzed far enough to yield
the estimate (1.4).
(This analysis is an application of the ``Poisson paradigm'' described by Alon and Spencer [A].)

A set of $k$ consecutive bit positions will be called a {\it $k$-block}.
There are $n-k+1$ distinct $k$-blocks.
A $k$-block will be said to be {\it active\/} if its rightmost position generates a carry and each of the remaining $k-1$ positions propagates a carry.
The event ``$\bfC_n \ge k$'' is clearly equivalent to the event ``there is at least one active $k$-block'',
which we shall denote $E_{n,k}$.
To estimate $\Pr[E_{n,k}]$, we shall use the following principles.
\medskip
\itemitem{(A-1)}
The probability that a given $k$-block is active is $1/2^{k+1}$.
\itemitem{(A-2)}
If a set of $k$-blocks includes two that overlap, then they cannot all be active.
If no two overlap, then they are independent.
\medskip

We shall show that (1.4) follows from these two principles.
Let
$$k_1 = \lceil 2\log_2 n\rceil.$$
For $k > k_1$, we have $\Pr[E_{n,k}] \le (n-k+1)/2^{k+1} = O(1/n)$ by (A-1) and Markov's inequality.
We also have $1 - e^{n/2^{k+1}} = O(1/n)$ by the power series $e^x = 1 + O(x)$, valid
for $x\to 0$.
Thus we have (1.4) for $k>k_1$.

For $k\le k_1$, we shall estimate $\Pr[E_{n,k}]$ using inclusion-exclusion, using (A-1) and (A-2).
We have
$$\eqalignno{
\Pr[E_{n,k}] &= \sum_{j\ge 1} {n - j(k-1)\choose j} {(-1)^{j-1} \over 2^{(k+1)j}} \cr
& = 1 - \sum_{j\ge 0} {n - j(k-1)\choose j} {(-1)^{j} \over 2^{(k+1)j}}, &(2.1)\cr
}$$
since there are just ${n - j(k-1)\choose j}$ ways to choose $j$ non-overlapping $k$-blocks
in the $n$ bit-positions.
Let 
$$k_0 = \left\lfloor \log_2 \left({3n \over 2 \log n - 6\log\log n}\right)\right\rfloor,$$
so that 
$${1\over 3}\log n - \log\log n \le  {n\over 2^{k_0+1}} \le {2\over 3}\log n - 2\log\log n,$$
$$e^{-n/2^{k_0+1}} = O\left({\log n\over n^{1/3}}\right)$$ 
and
$$e^{n/2^{k_0+1}} = O\left({n^{2/3} \over (\log n)^2}\right).$$ 
We shall begin by assuming $k\ge k_0$ (as well as $k\le k_1$).
Let 
$$j_0 = \left\lceil {(2e^2/3)} \log n\right\rceil.$$
We shall break the sum in (2.1) at $j_0$:
$$\Pr[E_{n,k}] = 1 - \sum_{0\le j \le j_0} {n - j(k-1)\choose j} {(-1)^{j} \over 2^{(k+1)j}}
 - \sum_{j>j_0} {n - j(k-1)\choose j} {(-1)^{j} \over 2^{(k+1)j}}. \eqno(2.2)$$
 We bound the magnitude of the second sum in (2.2) by using
 ${n - j(k-1)\choose j} \le {n\choose j} \le (en/j)^j$,
which yields 
$$\eqalignno{
\left\vert\sum_{j>j_0} {n - j(k-1)\choose j} {(-1)^{j} \over 2^{(k+1)j}}\right\vert
&\le \sum_{j>j_0} \left({en\over j2^{k+1}}\right)^j \cr
&\le \sum_{j>j_0} \left({1\over e}\right)^j \cr
& = O\left({1\over n^{2e^2/3}}\right). &(2.3)\cr
}$$

For the first sum in (2.2), we estimate the binomial coefficient by 
${n - j(k-1)\choose j} = (n^j/j!)\(1 + O(jk/n)\)^j = 
(n^j/j!)\(1 + O\((\log n)^3 / n\)\)$:
$$\sum_{0\le j \le j_0} {n - j(k-1)\choose j} {(-1)^{j} \over 2^{(k+1)j}}
= \sum_{0\le j \le j_0} {1\over j!} \left( {-n \over 2^{(k+1)}}\right)^j
\left(1 + O\left((\log n)^3\over n\right)\right).$$
The presence of the $O$-term in the summand prevents us from exploiting cancellation,
so to obtain an error bound for the sum we consider the magnitudes of the summands:
$$\eqalign{
\sum_{0\le j \le j_0} {n - j(k-1)\choose j} {(-1)^{j} \over 2^{(k+1)j}}
&= \left(\sum_{0\le j \le j_0} {1\over j!} \left( {-n \over 2^{(k+1)}}\right)^j\right)
 + O\left({(\log n)^3\over n} \sum_{0\le j \le j_0} {1\over j!} \left( {n \over 2^{(k+1)}}\right)^j\right) \cr
&= \left(\sum_{0\le j \le j_0} {1\over j!} \left( {-n \over 2^{(k+1)}}\right)^j\right)
 + O\left({(\log n)^3\over n} \, e^{n/2^{k+1}}\right) \cr
&= \left(\sum_{0\le j \le j_0} {1\over j!} \left( {-n \over 2^{(k+1)}}\right)^j\right)
 + O\left({\log n\over n^{1/3}} \right). \cr
}$$
Extending the sum from $j\le j_0$ to $j<\infty$ yields
$$\sum_{0\le j \le j_0} {n - j(k-1)\choose j} {(-1)^{j} \over 2^{(k+1)j}}
= \left(e^{-n/2^{k+1}} - \sum_{j>j_0} {1\over j!} \left( {-n \over 2^{(k+1)}}\right)^j\right)
 + O\left({\log n\over n^{1/3}} \right).$$
We bound the magnitude of this sum just as we did that of the second sum in (2.2), to obtain
$$\eqalignno{
\sum_{0\le j \le j_0} {n - j(k-1)\choose j} {(-1)^{j} \over 2^{(k+1)j}}
&= e^{-n/2^{k+1}}  + O\left({1\over n^{2e^2/3}}\right)
 + O\left({\log n\over n^{1/3}} \right) \cr
&= e^{-n/2^{k+1}}  + O\left({\log n\over n^{1/3}} \right). &(2.4)\cr
}$$
Substituting (2.3) and (2.4) in (2.2), we obtain (1.4) for $k_0 \le k \le k_1$.

Finally, we consider $k < k_0$.
We use the fact that $\Pr[E_{n,k}]$ is a non-increasing function of $k$, so that
$$1 \ge \Pr[E_{n,k}] \ge \Pr[E_{n,k_0}] = 1 - e^{-n/2^{k_0+1}}  + O\left({\log n\over n^{1/3}} \right)
= 1  + O\left({\log n\over n^{1/3}} \right).$$
This yields (1.4) for the remaining values of $k$.
\sk

\heading{3. The Algorithm for Multiplication}

In this section we shall describe in more detail the algorithm presented in the Introduction.
It will be most convenient to describe these algorithms in the language of hardware, implemented
as circuits built from gates interconnected by wires, but this is of course equivalent to
a description in the language of software for a parallel computer, such as that used by Claus [C]
and Knuth [K].

To begin, let us assume that $M$ is given its unique conventional binary representation
$M = \sum_{0\le j\le d} m_j \, 2^j$,
in which all digits $m_0 = 1, m_1, \ldots, m_d = 1$ are either $0$ or $1$.
As before, let $s_1 = 0 < s_2 < \cdots < s_c = d$ denote the positions of the $1$s.
Our first step will be to specify the encodings of the numbers $W_1, W_2, \ldots, W_c$
as sequences of bits.
The input $V = \sum_{0\le l\le n-1} v_l \, 2^l$ will be received using $n$ bits $v_0, v_1, \ldots, v_{n-1}$
as usual.
Since $V$ is an $n$-bit number (in the range $[0,2^n - 1]$) and $M$ is a $(d+1)$-bit number (in the range $[0, 2^{d+1} - 1]$), their product $Z = M\cdot V$ is an $(n+d+1)$-bit number
(in the range $\big[0, (2^n - 1)(2^{d+1} - 1)\big] \subseteq [0, 2^{n+d+1} - 1]$).
Thus it will suffice to represent all numbers produced during the execution of the algorithms 
(the output $Z$ and all intermediate results) using $n+d+1$ bits, and to perform all additions
(both carry-save and carry-propagate) modulo $2^{n+d+1}$.
Thus we shall represent each $W_i$ (for $1\le i\le c$) by the $n+d+1$ bits in its conventional
binary representation: $W_i = \sum_{0\le l\le n+d} w_{i,l} \, 2^l$.
Since $W_i = 2^{s_i} \, V$, we have $w_{i,l} = v_{l-s_i}$ if $s_i \le l\le n-1+s_i$, and
$w_{i,l} = 0$ if $0\le l\le s_i - 1$ or $n+s_i \le l\le n+d$.

In the first stage of the algorithm, we reduce the $c$ summands $W_1, W_2, \ldots, W_c$
to two summand $X$ and $Y$ by means of carry-save adders.
Each {\it carry-save adder\/} consists of $n+d+1$ ``full adders'', one for each position
in the numbers being added.
A {\it full adder\/} is a pair of gates that takes three input signals (say $f$, $g$ and $h$)
and produces two output signals.
The {\it sum\/} output is the parity (that is, the sum $f\oplus g\oplus h$ modulo $2$) of 
the three inputs.
The {\it carry\/} output is the majority ($(f\land g)\lor (f\land h)\lor (g\land h)$) of the three inputs.
The parity and majority are symmetric functions of the three inputs, so when specifying what signals
should be fed into a full adder, we do not need to specify which signal goes into which input.
The $n+d+1$ full adders in a carry-save adder reduce three summands
(say $F = \sum_{0\le l\le n+d} f_l \, 2^l$, $G = \sum_{0\le l\le n+d} g_l \, 2^l$ and 
$H = \sum_{0\le l\le n+d} h_l \, 2^l$) to two summands 
(say $A = \sum_{0\le l\le n+d} a_l \, 2^l$
and $B = \sum_{0\le l\le n+d} b_l \, 2^l$) as follows.
The signals $f_l$, $g_l$ and $h_l$ are fed into the inputs of the full adder in position $l$
(for $0\le l\le n+d$).
The sum outputs of the full adders become the bits of the summand $A$:
$a_l = \parity(f_l,g_l,h_l)$ for $0\le l\le n+d$.
Finally, the carry outputs of the full adders become, after being shifted left one position,
the bits of the summand $B$: $b_l{l+1}= \majority(f_l,g_l,h_l)$ for $0\le l\le n+d-1$
(the carry output from the full adder in the leftmost position is ignored) and $b_0 = 0$
(a $0$ bit is shifted into the rightmost position of $B$).

After the $c$ summands $W_1, W_2, \ldots, W_c$ have been reduced to two summands
$X$ and $Y$ by $c-2$ full adders in the first stage, the summands $X$ any $Y$ are added
by a carry-propagate adder in the second stage.
Like a carry-save adder, a carry-propagate adder can be built from $n+d+1$ full adders,
one for each position in the numbers being added.
Two of the inputs of the full adder in position $l$ (for $0\le l\le n+d$) are provided by the 
appropriate bits $x_l$ and $y_l$ of the numbers $X = \sum_{0\le l\le n+d} x_l \, 2^l$ and
$Y = \sum_{0\le l\le n+d} y_l \, 2^l$.
But in this case the third input of the  full adder in position $l$ is fed from the 
carry output of the full adder in position $l-1$ for $1\le l\le n+d$, and is fed the constant $0$
for $l=0$ (the carry output from the full adder in position $n+d$ is ignored).
The $n+d+1$ bits of the final product $Z$ are then produced at the sum outputs of the full adders.

This description of a carry-propagate adder gives an adequate picture of the production
of the outputs, but it is not convenient for the analysis of the longest carry propagation chain,
for which we must distinguish between between the {\it generation\/} of carries and their 
{\it propagation}, rather than merely their production.
To make the generation and propagation of carries more explicit, we will replace the full adders
in the second stage by ``half adders''.
A {\it half adder\/} is obtained from a full adder by substituting the constant $0$ for one of its
three inputs.
The resulting device consists of a pair of gates, one of which computes the sum output as the
parity (that is, the ``exclusive-OR'') of the two remaining inputs, and the other of which computes
the carry output as the conjunction (that is, the ``AND'') of the inputs.
If we replace each full adder in the second stage with a half adder, then the carry output
of each half adder will indicate whether a carry is generated at that position (that is, whether
both $x_l$ and $y_l$ are $1$s for that value of $l$), and the sum output will indicate whether a carry would be propagated by that position (that is, whether exactly one of $x_l$ and $y_l$ is a $1$).

We conclude this section with a discussion of how to adapt the algorithm presented above for a recoded multiplier, in which some of the digits $m_j$ ($0\le j\le d-1$) may be $\overline{1} = -1$.
(We note that since $M$ is positive, we must have $m_d = 1$.)
We could do this by changing the adder that incorporates the contribution $W_j$ into a subtracter,
but this might require analysis of subtraction, rather than addition, in the second stage (if 
$m_0 = \overline{1}$).
It will be more convenient to preserve the structure of the adders in both stages, and to alter the definition of the $W_j$ ($0\le j\le d$) to be a positive integers whose sum is congruent to $M\cdot V$ 
modulo $2^{n+d}$.
To do this, we make the following changes for each $j$ ($0\le j\le d-1$) such that $m_j = \overline{1}$:
\medskip
\itemitem{(a)}
Change the leftmost $d-s_j$ bits $w_{j,n+d-1}, \ldots, w_{j,n+s_j}$ of $W_j$ from $0$s to $1$s.
\itemitem{(b)}
Change the $n$ bits $w_{j,n+s_j-1}, \ldots, w_{j,s_j}$ of $W_j$ from $v_{n-1},\ldots, v_0$ to their complements $\lnot{v}_{n-1},\ldots, \lnot{v}_0$.
\itemitem{(c)}
Change the bit $w_{j+1,s_j}$ of $W_{j+1,s_j}$ from $0$ to $1$.
\medskip
\noindent These changes result in a contribution congruent to $m_j \, 2^{s_j} \, V$ modulo $2^{n+d}$, since changes (a) and (b) form $2^{s_j}$ times the $1$s complement of $V$, and change (c) converts
the $1$s complement to the $2$'s complement.
\sk

\heading{4. The Analysis of Multiplication}

We begin by deriving the principles, analogous to (A-1) and (A-2), that will allow us to analyze multiplication.
A {\it $k$-block\/} is a sequence of contiguous bit positions among the $n+d$ positions of numbers
modulo $2^{n+d}$.
Thus there are just $n+d-k-1$ distinct $k$-blocks, with the rightmost position of the rightmost $k$-block being position $0$, and the leftmost position of the leftmost $k$-block being position $n+d-k$.
We shall say that a $k$-block is {\it active\/} if, in the final addition in the second stage, its rightmost position generates a carry and its remaining $k-1$ positions propagate a carry.
Whether or not a $k$-block is active depends on the input bits not only in its $k$ positions,
but also in up to $d$ positions to its right.
These $d$ or fewer positions will be called the {\it extension\/} of the $k$-block, and the $k$-block together with its extension will be called an {\t extended\/} $k$-block.
(The $d$ rightmost $k$-blocks will have fewer than $d$ positions in their extensions, since there are
fewer than $d$ positions to their right.)

The inputs to the final addition are computed by circuits composed of three-input parity and majority gates, one-input inverters, and zero-input constant gates.
Furthermore, constant gates occur only in the circuits computing the rightmost $d$ and leftmost $d$ positions (positions $0$ through $d-1$ and positions $n$ through $n+d-1$).
A $k$-block will be called {\it marginal\/} if it r its extension overlap the rightmost or leftmost $d$ positions.
Thus there are $3d$ marginal $k$-blocks.
A $k$-block will be called {\it central\/} if it is not marginal.
\medskip
\itemitem{(M-1)}
The probability that a central $k$-block is active is $1/2^{k+1}$.
\medskip

Suppose the rightmost position of the $k$-block is position $l$ ($2d\le l\le n-d-k$).
For the rightmost position to generate a carry, the values of both $x_l$ and $y_l$ must be $1$.
The value of $y_l$ depends on the inputs $v_{l-1}, \ldots, v_{l-d}$, and it is computed from them by a circuit composed of three-input parity and majority gates and one-input inverters.
These gates compute {\it self-dual\/} Boolean functions: if the arguments of a self-dual function are complemented, then the value of the function is also complemented.
The class of self-dual functions is closed under composition, so $y_l$ is a self-dual function
of the inputs $v_{l-1}, \ldots, v_{l-d}$.
If the arguments of a self-dual function are independent unbiassed bits, then the value of the function
is also an unbiassed bit.
Thus the probability that $y_l = 1$ is $1/2$.
The value of $x_l$ depends on the input $v_l$ as well as the $d$ inputs to its right, and we have
$$x_l = v_l \oplus \ph(v_{l-1}, \ldots, v_{l-d}),$$
where $\ph$ is some $d$-adic Boolean function.
Since $v_l$ is an unbiassed bit independent of $v_{l-1}, \ldots, v_{l-d}$, $x_l$ is an unbiassed bit independent of $y_l$.
Thus the probability that position $l$ generates a carry is $1/4$.

For each of the remaining $k-1$ positions of the $k$-block to propagate a carry, we must have
$x_{l+1}\oplus y_{l+j} = 1$ for $1\le j\le k-1$.
As between $x_{l+1}$ and $y_{l+j} = 1$, only $x_{l+1}$ depends on $v_{l+j}$ and, as above, we have
$$x_{l+j} = v_{l+j} \oplus \ph(v_{l+j-1}, \ldots, v_{l+j-d}).$$
Thus each $x_j$ is an unbiassed bit independent of the bits to its right, so the probability of each 
of the remaining $k-1$ bits propagating a carry is $1/2^{k-1}$, and the probability that a central $k$-block is active is $(1/4)(1/2^{k-1}) = 1/2^{k+1}$.
\medskip
\itemitem{(M-2)}
The probability that a marginal $k$-block is active is at most $2^d/2^{k}$.
\medskip

The analysis of (M-1) applies to  the $k-d$ or more positions of the $k$-block that do not overlap
the rightmost $2d$ or leftmost $d$ positions.

We shall say that two $k$-blocks are {\it strongly non-overlapping\/} if they, together with their extensions,
are non-overlapping, and that they are {\it weakly overlapping\/} if they are non-overlapping, but one overlaps the extension of the other.
\medskip
\itemitem{(M-3)}
If two $k$-blocks are overlapping, they cannot both be active.
\medskip

This holds because at each position, generating a carry and propagating a carry are exclusive events.
\medskip
\itemitem{(M-4)}
If a $k$-block $B$ lies to the right of, and is strongly non-overlapping,  a $k$-block $A$ then the event that $B$ is active is independent of the event that $A$ is active.
\medskip

This holds because the activities of strongly non-overlapping $k$-blocks depend on disjoint sets of inputs.
\medskip
\itemitem{(M-5)}
If a $k$-block $B$ overlaps the extension of a $k$-block $A$, but does not overlap $A$ itself,
then the probability that $B$ is active, given that $A$ is active, is at most $2^d/2^{k+1}$.
\medskip

The analysis of (M-1) applies to the $k-d$ or more rightmost positions of $B$ that do not overlap
$A$ or its extension.

We shall show that (1.4) follows from these five principles.
As before, we let
$$k_1 = \lceil 2\log_2 n\rceil.$$
Then (1.4) follows for $k>k_1$, since from (M-1), (M-2) and Markov's inequality,
$\Pr(E_{n,k})$ is $O(1/n)$, as is $1 - e^{-n/2^{k+1}}$.

For $k\le k_1$, we again let
$$k_0 = \left\lfloor \log_2 \left({3n \over 2 \log n - 6\log\log n}\right)\right\rfloor,$$
and begin by assuming that $k\ge k_0$ (as well as $k\le k_1$).
Using (M-2) and Markov's inequality, the probability that any marginal $k$-block is active is at most
$3d2^d/2^k = O(\log n/n)$.
Thus we may ignore marginal $k$-blocks, and turn our attention to estimating the probability of the event
$E'_{n,k}$ that some central $k$-block is active.
For this, we shall again use inclusion-exclusion:
$$\Pr(E'_{n,k}) = \sum_{j\ge 1} \sum_{B_1,\ldots,B_j} (-1)^{j-1}\Pr(B_1,\ldots,B_j\hbox{\ all active\ }),
\eqno(4.1)$$
where the sum is over all lists $(B_1,\ldots,B_j)$ of $j$ central $k$-blocks, with $B_{i+1}$ to the right
of $B_i$ for $1\le i\le j-1$.
By (M-3), we may also assume that $B_1,\ldots,B_j$ are pairwise non-overlapping.

We shall partition the contributions to the double sum in (4.1) into two parts,
$$\Pr(E'_{n,k}) = \Si_\rmi + \Si_\rmii,$$
where $Si_\rmi$ denotes the sum of the contributions from lists $B_1,\ldots,B_j$ that are pairwise strongly non-overlapping, and $\Si_\rmii$ denotes the sum of the contributions from lists 
$B_1,\ldots,B_j$ for which at least one pair $B_i$, $B_{i+1}$ of successive $k$-blocks is weakly
overlapping.
The contributions to $\Si_\rmi$ will be completely analogous to those in the analysis of addition.
For the contributions to $\Si_\rmii$, we shall need to analyze the effects of weak overlaps,
but in this case it will suffice to consider only the magnitudes of the contributions, without making
any attempt to exploit cancellations.

For $\Si_\rmi$, the only difference from the analysis of addition is that now the extended $k$-blocks each have length $k+d$, and the number of positions into which $j$ of them must fit is now $n-2d$.
Thus the binomial coefficient that counts the number of ways that $j$ strongly non-overlapping central 
$k$-blocks can be chosen is ${(n-2d) - j(k+d-1)\choose j}$.
Since this quantity still satisfies the estimates
$${(n-2d) - j(k+d-1)\choose j} = {n^j\over j!}\left(1 + O\left({(\log n)^3\over n}\right)\right)$$
for $j\le j_0$, where again
$$j_0 = \left\lceil {(2e^2/3)} \log n\right\rceil,$$
and
$${(n-2d) - j(k+d-1)\choose j} \le {n^j\over j!} \le \left({en\over j}\right)^j$$
for all $j$, we can use (M-1) in the analysis of Section 2 to show that
$$\Si_\rmi = 1 - e^{-n/2^{k+1}} + O\left({\log n\over n^{1/3}}\right).$$

Turning to $\Si_\rmii$, we abandon any attempt to exploit cancellation among the terms, and merely sum bounds on their magnitudes.
We have
$$\abs{\Si_\rmii} \le \sum_{l\ge 1} \sum_{j\ge l+1} {j-1\choose l}
\sum_{1\le f_1,\ldots,f_l\le d} {(n-2d) - (j-l)(k+d-1)+g\choose j-l}
\left({1\over 2^{k+1}}\right)^{j-l} \left({2^d\over 2^{k+1}}\right)^{l},$$
where $g = f_1 + \cdots + f_l$.
Here $l$ denotes he number of values of $i$ ($1\le i\le j-1$) such that $B_{i+1}$ overlaps the extension of $B_i$, the binomial coefficient ${j-1\choose l}$ counts the number of ways in which these values of $i$ may be chosen, the parameters $f_1, \ldots, f_l$ denote the amounts of overlap, the binomial coefficient  ${(n-2d) - (j-l)(k+d-1)+g\choose j-l}$ counts the number of ways in which the $j-l$ $k$-blocks
or weakly overlapping sequences of $k$-blocks may be chosen, the factor $(1/2^{k+1})^{j-l}$ denotes the probability, following (M-1) and (M-4), that the $j-l$ $k$-blocks that do not overlap the extension
of a $k$-block to their left are all active, and the factor $(2^d/2^{k+1})^{l}$ bounds the probability, following (M-5), that the remaining $l$ $k$-blocks are all active.
Since the innermost sum has at most $d^l$ terms, each with $g\le ld$, the innermost binomial coefficient
is at most ${n\choose j-l}\le n^{j-l}/(j-l)!$ and we obtain
$$\eqalign{
\abs{\Si_\rmii} &\le \sum_{l\ge 1} \left({d2^d\over 2^{k+1}}\right)^l \sum_{j\ge l+1} {j-1\choose l}
{1\over (j-l)!} \left({n\over 2^{k+1}}\right)^{j-l} \cr
&= \sum_{l\ge 1} \left({d2^d\over 2^{k+1}}\right)^l \sum_{m\ge 1} {m+l-1\choose l}
{1\over m!} \left({n\over 2^{k+1}}\right)^{m}, \cr
}$$
where we have made the substitution $m = j-l$.
We shall show below that 
$$\sum_{m\ge 1} {m+l-1\choose l}
{x^m\over m!} \le (4x)^l e^x \eqno(4.2)$$
for $x\ge 1$ and $l\ge 1$.
Since $n/2^{k+1} \ge (1/3)\log n - \log\log n\ge 1$ and 
$e^{n/2^{k+1}} = O\(n^{2/3}/(\log n)^2\)$ for all sufficiently large $n$ and $k\ge k_0$, 
we obtain
$$\eqalign{
\abs{\Si_\rmii} &\le \sum_{l\ge 1} \left({d2^d\over 2^{k+1}}\right)^l 
\left({4n\over 2^{k+1}}\right)^l  e^{n/2^{k+1}}\cr
&=O\left({1\over n^{1/3}}\right). \cr
}$$
It remains to prove (4.2).
We have
$$\eqalign{
\sum_{m\ge 1} {m+l-1\choose l} {x^m\over m!}
&= {x\over l!}{d^l\over dx^l} \sum_{m\ge 1} {x^{m+l-1}\over m!} \cr
&= {x\over l!}{d^l\over dx^l} \,\, x^{l-1}\sum_{m\ge 1} {x^{m}\over m!} \cr
&= {x\over l!}{d^l\over dx^l} \,\, x^{l-1}(e^x - 1) \cr
&= {x\over l!}\sum_{0\le s\le l}{l\choose s}\left({d^s\over dx^s} \,\, x^{l-1}\right)
\left({d^{l-s}\over dx^{l-s}}(e^x - 1) \right) \cr
&= {x\over l!}\sum_{0\le s\le l-1}{l\choose s} {l-1\choose s}s! \,\, x^{l-1-s}
e^x. \cr
}$$
Since $x\ge 1$, we have $x^{l-1-s} \le x^{l-1}$. 
Using the further inequalities
$s!\le l!$ and
$$\eqalign{
\sum_{0\le s\le l-1} {l\choose s} {l-1\choose s}
&\le \sum_{0\le s\le l} {l\choose s} {l\choose s} \cr
&= {2l\choose l} \cr
&\le 4^l, \cr
}$$
we obtain (4.2).
This competes the proof of (1.4) for $k_0\le k\le k_1$.

Finally, we must consider $k<k_0$.
Again as in the analysis of addition, the fact that $\Pr[E'_{n,k}]$ is a non-increasing function of $k$,
together with the bound (1.4) for $k=k_0$ yields (1.4) for the remaining values of $k$.
\sk

\heading{5. Conclusion}

In this paper we have shown that the distribution of the length of the longest carry propagation chain can be analyzed using what Alon and Spencer [A] have called the ``Poisson paradigm''.
We have also show that this method of analysis can be used to show that a particular algorithm
for multiplication of a random integer by a fixed constant has, to within terms tending to zero as 
$n\to\infty$, the same distribution for the length of the longest carry chain in the final addition.
This algorithm is characterized by shifting over zeros in the multiplier, and by the use of a carry-save
adder to incorporate the contributions for all but the last two non-zero digits of the multiplier.
We should point out that our analysis does not appear to be applicable to either of two natural variants of this algorithm: one in which zeros are not shifted over, but cause a contribution of zero to be added using a carry-save adder (for in this case we cannot appeal to self-duality in the computation of the 
final summands), and one in which a carry-propagate adder is used for all additions (in which case it does not matter whether or not zeros are shifted over, for in this case the outputs of each adder depend
on an unbounded number of input bits to their right).
It remains an open question whether the result of this paper applies to either or both of these variants.
An apparently even more challenging problem is to determine whether or not the result of this paper
applies to the algorithm considered here when the multiplier is not a fixed integer, but is rather a 
random integer with the same distribution as, but independent of, the multiplicand.
\sk

\heading{6. Acknowledgment}

The research reported here was supported in part by NSF Grant CCF 0430656.
\sk

\heading{7.  References}

\refbook A; N. Alon and J. H. Spencer;
The Probabilistic Method; (Second Ed.)
John Wiley \& Sons, 2000. 

\refinbook B; A. W. Burks, H. H. Goldstein and J. von Neumann;
``Preliminary Discussion of the Logical Design of an Electronic Computing Instrument'';
in: A.~H. Taub (Ed.); Collected Works of John von Neumann;
Macmillan, 1963, v.~5, pp.~34--79. 

\ref C; V. Claus;
``Die mittlere Additionsdauer eines Paralleladdierwerks'';
Acta Informatica; 2 (1973) 283--291.

\refbook I; A. Izsak;
Special Cases of Carry Propagation;
B.~S. Thesis, Department of Mathematics, Harvey Mudd College, 2007.

\ref K; D. E. Knuth;
``The Average Time for Carry Propagation'';
Nederl.\ Akad.\ Wettensch.\ Indag.\ Math.; 40 (1978) 238--242
 (reprinted in D.~E. Knuth, {\it Selected Papers on Analysis of Algorithms}, 
Center for the Study of Language and Information, Stanford University, 2000). 

\ref L1; M. Lehman;
``High-Speed Digital Multiplication'';
IRE Trans.\ Electronic Computers; 6 (1957) 204--205.

\ref L2; M. Lehman;
``Short-Cut Multiplication and Division in Automatic Binary Digital Computers'';
Proc. IEE; 105 B (1958) 496--504.

\ref P; N. Pippenger;
``Analysis of Carry Propagation in Addition: An Elementary Approach'';
J. Algorithms; 42 (2002) 317--333.

\ref R; G. W. Reitwiesner;
``Binary Arithmetic'';
Advances in Computers; 1 (1960) 232--308.

\ref T; K. D. Tocher;
``Techniques of Multiplication and Division for                                 
Automatic Binary Computers'';
Quart.\ J. Mech.\ Appl.\ Math.; 11 (1958) 364--384.

\ref W; C. S. Wallace;
``A Suggestion for a Fast Multiplier'';
IEEE Trans.\ Computers; 13:2 (1964) 14--17.

\bye